\documentclass[10pt]{article}

\usepackage{algorithm2e}
\usepackage{xcolor}

\SetAlFnt{\footnotesize}

\SetAlCapSty{xAlCapSty}

\SetCommentSty{xCommentSty}

\SetNlSty{mynlfont}{}{} 

\LinesNumbered

\SetSideCommentRight

\DontPrintSemicolon

\RestyleAlgo{algoruled}

\usepackage[thinlines]{easytable}

\usepackage{color}      
\usepackage{placeins}
\usepackage{indentfirst}
\usepackage{multirow}

\usepackage{amssymb}
\usepackage{amsmath,latexsym}
\usepackage{url}
\usepackage{enumerate}
\usepackage{graphicx}

\usepackage{booktabs}

\usepackage{subcaption}
\usepackage{cleveref}
\usepackage{mathtools}

\begin{document}

\title{\bf On Pitfalls in Accuracy Verification Using Time-Dependent Problems}

\author{Hiroaki Nishikawa}
\date{\today}
\maketitle

\begin{abstract}
In this short note, we discuss the circumstances that can
lead to a failure to observe the design order of discretization error convergence
in accuracy verification when solving a time-dependent problem. In particular,
we discuss the problem of failing to observe the design order of spatial
accuracy with an extremely small time step. The same problem is encountered even if the time step is reduced with grid refinement. 
These can cause a serious problem
because then one would wind up trying to find a coding error that does
not exist. This short note clarifies the mechanism causing this failure and provides a guide
for avoiding such pitfalls. \end{abstract}

\section{Introduction}
\label{intro}

Accuracy verification is an important process of checking a correct implementation of a numerical algorithm in a code \cite{Roache_verification_book,ThomasDiskinRumsey2008,OberkampfRoy_2010,tmr_link}, and it has been widely employed in the development of computational fluid dynamics codes \cite{ThomasDiskinRumsey2008,Roy_etal_IJNMF2007,aiaa_guide_vv:1998}. For example, one can verify the implementation of a second-order spatial discretization scheme in a code by directly measuring the order of discretization error (i.e., solution error) convergence over a series of consistently refined grids and demonstrating that it is indeed second order. The exact solution, which is required for computing discretization errors, can be created by the method of manufactured solutions \cite{OberkampfRoy_2010}. It is usually straightforward to perform such an accuracy verification, but as we report in this short note, accuracy verification needs to be performed carefully for time-dependent problems, or the design order of accuracy will not be observed. As we will show, we can easily fall into such a pitfall especially when trying to save time, for example, by measuring the discretization error after just one tiny time step. Intuitively, it seems to be a reasonable thing to do because then the time integration is virtually exact. However, this approach can result in one-order lower discretization error convergence. The problem is not necessarily resolved even if the time step is consistently reduced along with the grid refinement. These can cause a serious problem  because then one would waste time trying to find a coding error that does not exist. It is the objective of this short note to clarify the mechanism causing the failure to observe the design order of discretization error convergence in accuracy verification when solving time-dependent problems and provide a guide for avoiding such a pitfall.


\vspace{-0.3cm}
\section{A Second-Order Scheme}
\label{Fromm}

  \begin{figure}[t]
    \centering
          \begin{subfigure}[t]{0.32\textwidth}
        \includegraphics[width=\textwidth]{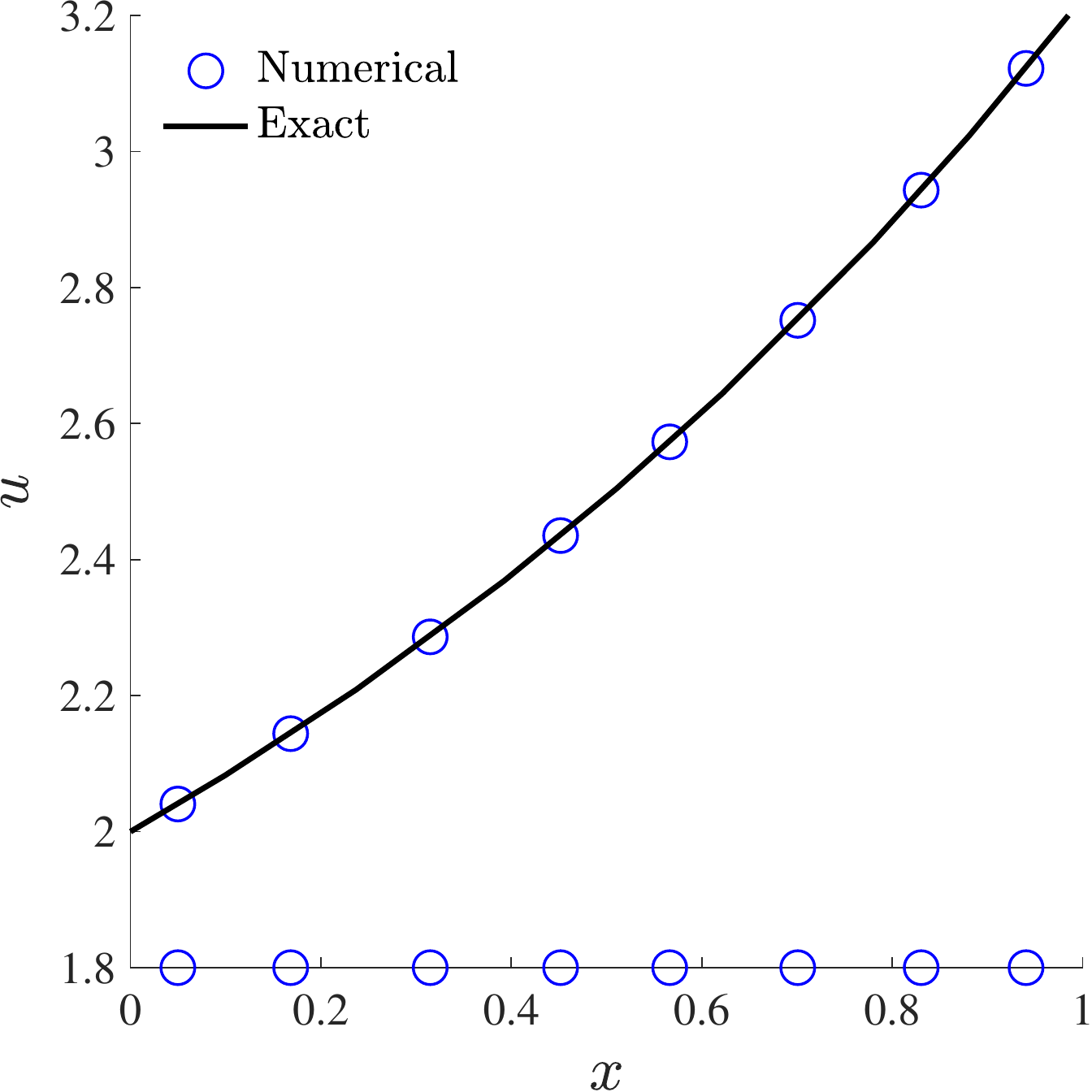}
          \caption{Solutions at $t=10^{-8}$ {\color{blue} after one time step}.}
          \label{fig:00_sol_irrg}
      \end{subfigure}
      \hfill
      \begin{subfigure}[t]{0.32\textwidth}
        \includegraphics[width=\textwidth]{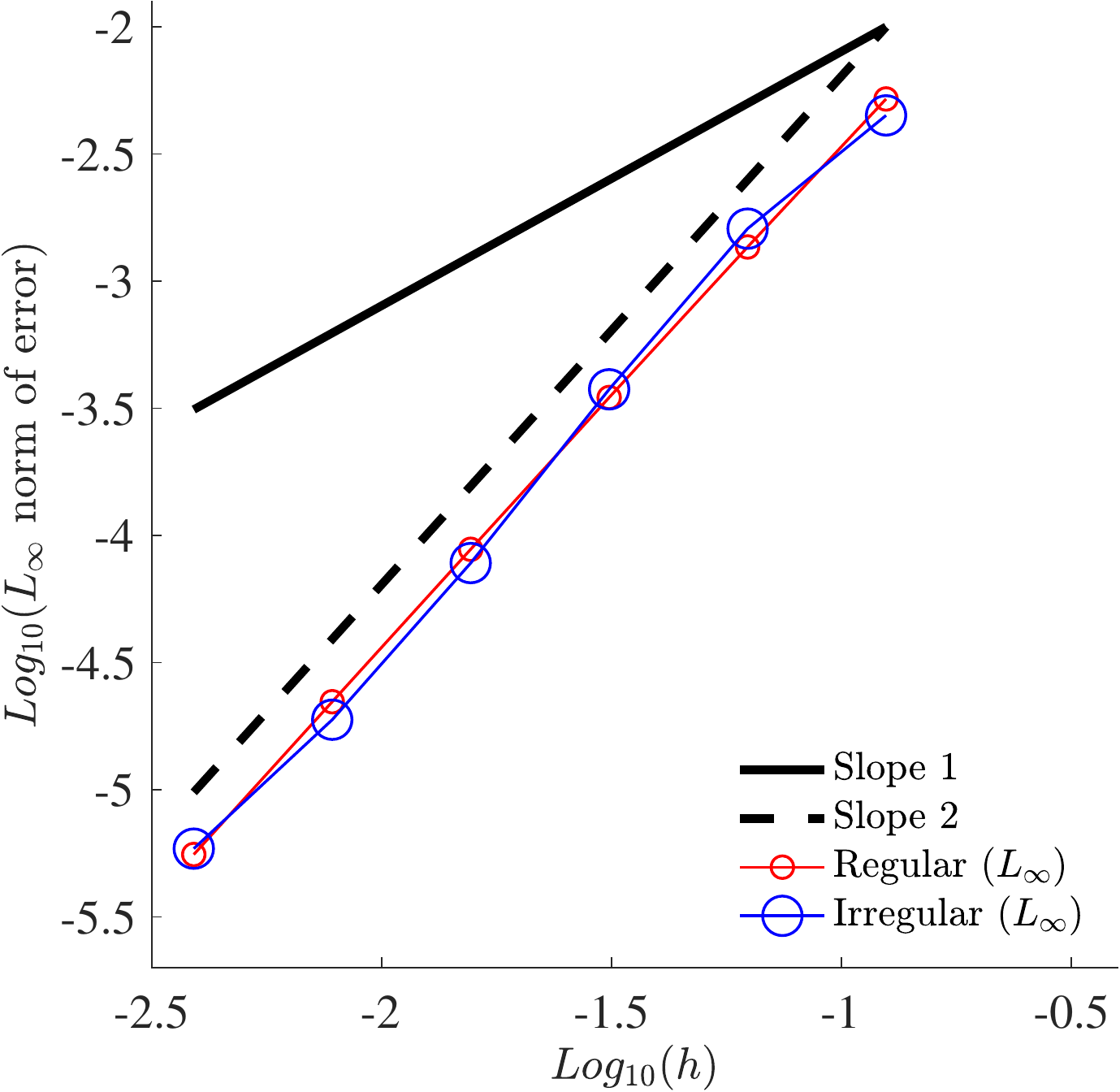}
          \caption{Steady problem.}
          \label{fig:0001_Li}
      \end{subfigure}
      \hfill
      \begin{subfigure}[t]{0.32\textwidth}
        \includegraphics[width=\textwidth]{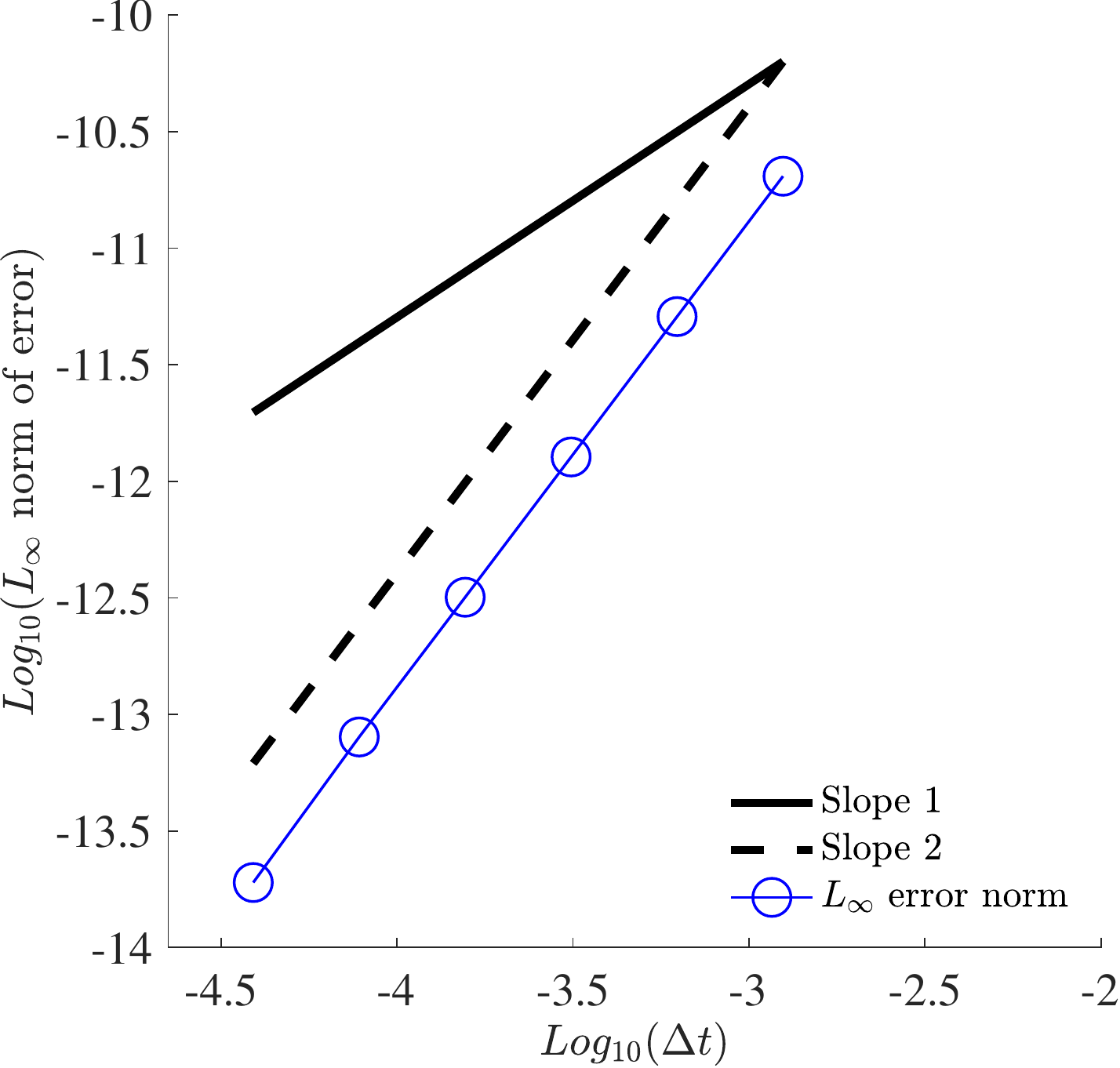}
          \caption{Time integration problem.}
          \label{fig:0002_Li}
      \end{subfigure} 
      \begin{subfigure}[t]{0.37\textwidth}
        \includegraphics[width=\textwidth]{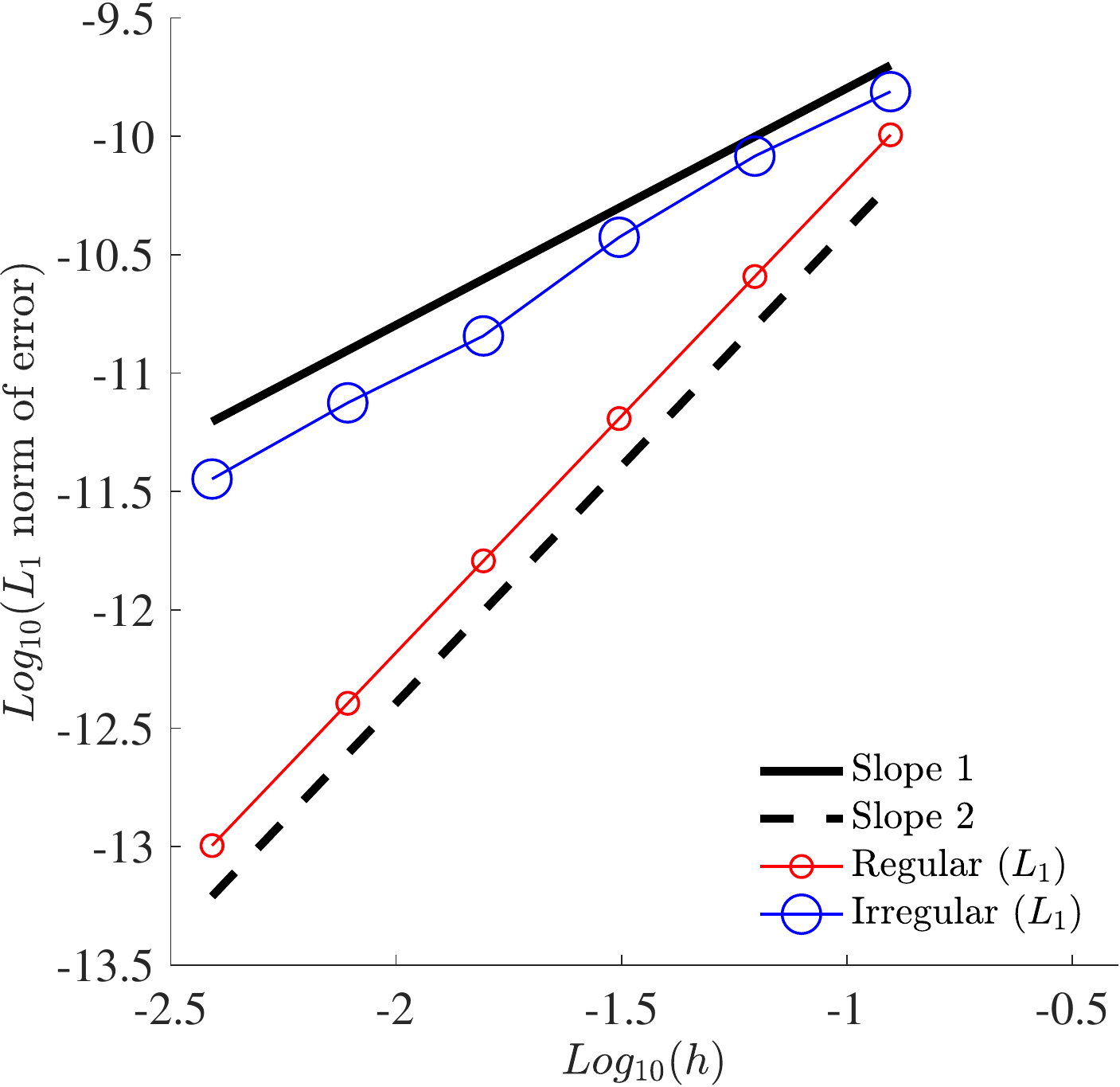}
          \caption{$L_1$ errors at $T_f=10^{-8}$.}
          \label{fig:00_L1}
      \end{subfigure}
      \begin{subfigure}[t]{0.37\textwidth}
        \includegraphics[width=\textwidth]{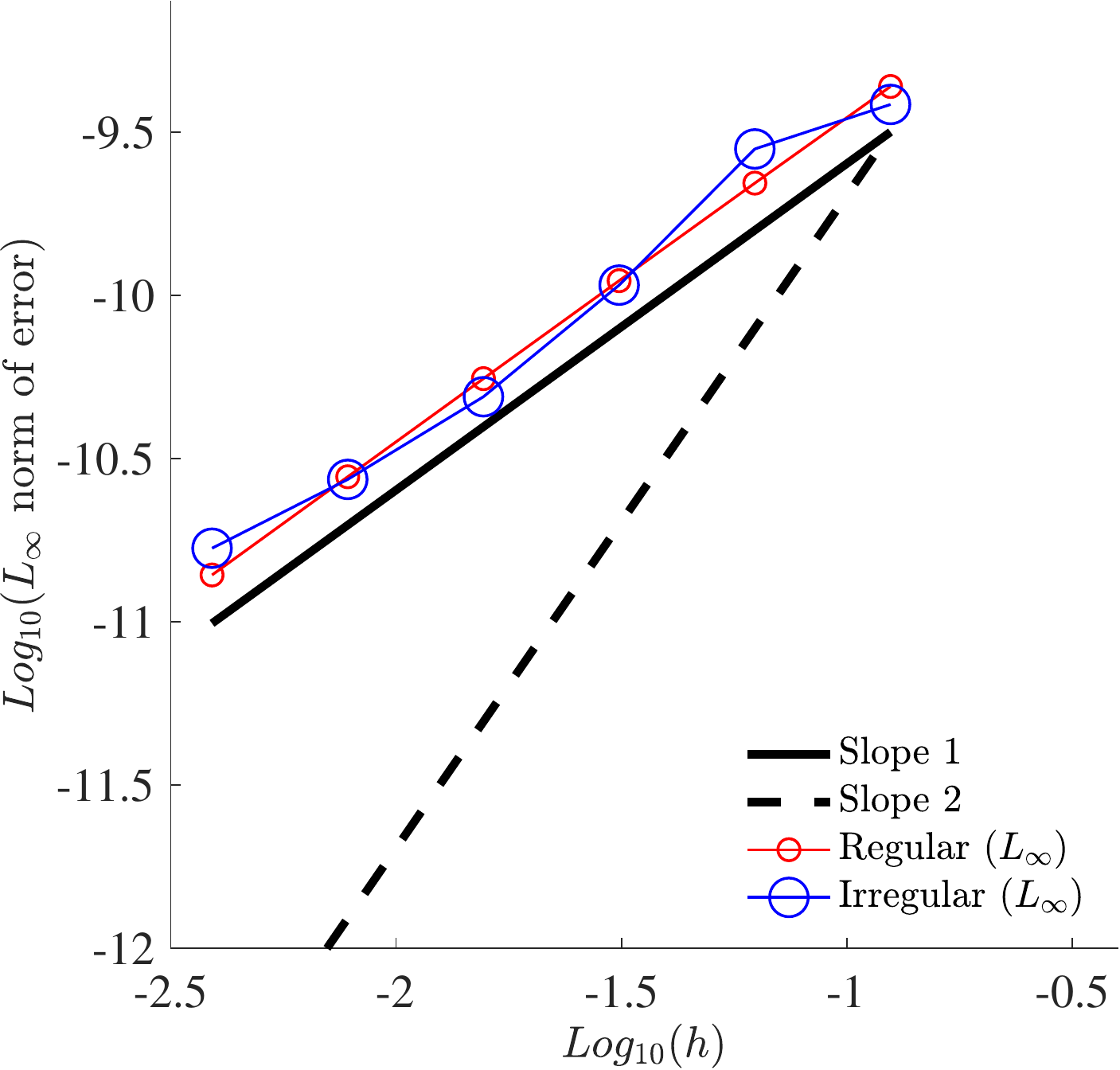}
          \caption{$L_\infty$ errors at $T_f=10^{-8}$.}
          \label{fig:00_Li}
      \end{subfigure} 
            \caption{
\label{fig00_results}%
Accuracy verification test for a time-dependent problem {\color{blue} with $a=1$} on regular and irregular grids of $N = 2^{k+2}$ cells, $k=1,2, \cdots, 6$: 
(a) numerical and exact solutions on the coarsest grid at $t = T_f = 10^{-8}$ {\color{blue} after one time step}, where circles on the $x$-axis indicate that the grid is not uniform, 
(b) results for a steady problem, (c) results for an ordinary differential equation, (d) $L_1$ errors for a time-dependent problem at $\Delta t = 10^{-8}$, and (e) $L_\infty$ errors for a time-dependent problem at $\Delta t = 10^{-8}$.
} 
\end{figure}

Consider an advection equation in one dimension, 
\begin{eqnarray}
\partial_t u + a \, \partial_x u = s(x,t),
 \,\,\, \mbox{in} \,\,\, x \in (0,1),
\end{eqnarray}
where $u$ is a solution variable, $a $ is a positive constant advection speed, $s(x,t)$ is a forcing term defined by $s(x,t) = \partial_t u^{exact} + a \, \partial_x u^{exact} 
$ for the following exact solution: 
\begin{eqnarray}
u^{exact}(x,t) = 1 + e^{0.8 x - 0.35 t},
\end{eqnarray}
which is plotted {\color{black} for $a=1$} in Figure \ref{fig:00_sol_irrg} at $t=10^{-8}$ on an irregular grid with 8 cells. {\color{black} Suppose we store a set of numerical solutions at cell centers, $x_j$, $j=1, 2, 3, \cdots, N$, in an arbitrarily-spaced grid of $N$ cells. Then, we have implemented a second-order accurate scheme using the second-order SSP-RK scheme \cite{SSP:SIAMReview2001} for time integration:
\begin{eqnarray}
u_j^{(1)} &=& u_j^{n}  -  \frac{\Delta t}{h_j} \left[   {Res}_j^n -  s( x_j, t^n ) h_j  \right] ,   \\ [2ex] 
u_j^{n+1} &=& \frac{ u_j^{n} + u_j^{(1)} }{2}   - \frac{\Delta t}{2 h_j} \left[     {Res}_j^{(1)}  - s( x_j, t^n + \Delta t )  h_j \right],
\end{eqnarray}
where $u_j^n$ is a numerical solution at $ t=t^n$ and at $x=x_j$, $h_j$ is the volume of the cell $j$, $n$ denotes the time level, $\Delta t$ is a time step, the spatial discretizations $Res_j^n$ and ${Res}_j^{(1)}$ are evaluated with $u_j^{n}$ at $t=t^n$ and $u_j^{(1)}$ at $t = t^n + \Delta t$, respectively.} The spatial discretization is defined by a second-order upwind scheme (Fromm's scheme \cite{Fromm:JCP1968}):
\begin{eqnarray} 
Res_j^n = f_{j+1/2} - f_{j-1/2}, \quad 
 f_{j+1/2} = a \left[  u_j^n + \left( \frac{ u_{j+1}^n - u_{j-1}^n}{ x_{j+1} - x_{j-1} }  \right) \frac{h_j}{2} \right], 
 \quad  
 f_{j-1/2} = a \left[   u_{j-1}^n + \left( \frac{ u_{j}^n - u_{j-2}^n}{ x_{j} - x_{j-2} }  \right) \frac{h_{j-1}}{2}  \right], 
 \label{original_scheme}
\end{eqnarray}
where the fluxes for the cells adjacent to the left and right boundaries are defined by
\begin{eqnarray}
 f_{1/2} =  a \left[  u^{exact}(0,t^n)  \right], \quad
 f_{3/2} =  a \left[  u_1^n + \left( \frac{ u_{2}^n - u_{1}^n}{ x_{2} - x_{1} }  \right) \frac{h_1}{2}  \right], \quad
 f_{N+1/2} = a \left[  u_N^n + \left( \frac{ u_{N}^n - u_{N-1}^n}{ x_{N} - x_{N-1} }  \right) \frac{h_N}{2}  \right].
\end{eqnarray}
It is well known that this spatial discretization is second-order accurate on irregular grids in one dimension because the algorithm is exact everywhere for linear solutions \cite{Boris_Jim_NIA2007-08}. One can easily verify second-order accuracy by solving a steady problem, ${Res}_j -  s( x_j, 0) h_j = 0$ over a series of consistently refined regular or irregular grids of  $N = 2^{k+2}$ cells, $k=1,2,3, \cdots, 6$, and measuring the error against the exact solution $u^{exact}(x,0)$. See Figure \ref{fig:00_sol_irrg} for the coarsest irregular grid. Our results obtained for $a=1$ are shown in Figure \ref{fig:0001_Li}. Also, it is straightforward to verify second-order time accuracy of the SSP-RK scheme by solving an ordinary differential equation $d u  / dt = s(x,t),$ with ${Res}_j = 0$ toward a final time $t = T_f$ over a set of consistently refined time steps and measuring the error against $u^{exact}(x,t)$ with initial values given by $u^{exact}(x_j,0)$. Our results are shown in Figure \ref{fig:0002_Li}, which were obtained with $T_f = 0.01 h_c $ and $\Delta t = 0.01 h$, where $h = 1/2^{k+2}$ and $h_c = 1/8$. Clearly, second-order error convergence is observed in both cases.

However, even though we have verified the spatial discretization and the time integration scheme independently, we encounter lower-order error convergence when we perform accuracy verification using a time-dependent problem. Consider, again, a simple case with $a=1$. First, we computed a solution at $T_f = 10^{-8}$ with a constant time step $\Delta t = 10^{-8}$, for the same regular and irregular grids used in the steady problem. Figure \ref{fig:00_sol_irrg} shows numerical and exact solutions on the coarsest irregular grid. Here, the time step is so small that the time integration is essentially exact after just one time step, and we would expect to observe second-order error convergence. However, we obtained second-order accuracy only in the $L_1$ norm: $\sum_{j=1}^N | u_j -u^{exact}(x_j,T_f) | / N$, and only on regular grids, as can be seen in Figure \ref{fig:00_L1}. {\color{black} Otherwise, we observe only first-order accuracy in $L_\infty \,\, \mbox{norm} = \max_{j = 1,2,3, \cdots, N}  | u_j -u^{exact}(x_j,T_f) |$ for the regular grids as in Figure \ref{fig:00_L1}, and in both $L_1$ and $L_\infty$ norms for the irregular grids, as in Figure \ref{fig:00_Li}.} In the case of regular grids, the maximum error is found at a cell adjacent to a boundary. Note that the errors are extremely small and thus the solutions are very accurate; but the rate of error convergence {\color{black} is not the expected rate}. 

We have also run this case for a longer time up to $T_f = 0.01 h_c$ with a larger consistently refined time step $\Delta t = 0.01 h$, $h = 1/2^{k+2}$, $k=1,2, \cdots, 6$. In this case, we used $ T_f/ \Delta t = h_c/ h = 2^{k-1} = 1, 2, 4, 8, 16, 32$ time steps for the grids of $8, 16, 32, 64, 128, 256$ cells, respectively. But the results are essentially the same as those shown in Figures \ref{fig:00_L1} and \ref{fig:00_Li} with larger errors (therefore not shown). We have encountered similar results in accuracy verification studies for second- and third-order schemes in three dimensions, where first- and second-order errors were observed, respectively, when applied to a time-dependent problem with an extremely small time step over a single time step $T_f = \Delta t$. 
o perform the verification over just one or a few small time steps.

We emphasize that the above scheme is very well known to be second-order
accurate on irregular grids in one dimension. But these first-order
accurate results should not be taken as a sign of a coding error, or one would be chasing an error that does not exist. As we will
show in the next section, these problems are originated from a one-order
lower truncation error on irregular grids, which include stencils adjacent
to a boundary in a regular grid, and, in some way, are a manifestation
of coarse grid behavior. Later, we will show that there is a simple approach
that allows us to avoid these pitfalls.

\vspace{-0.3cm}
\section{Causes of Lower-Order Behaviors}

To reveal the mechanism of the lower-order behavior, we focus on cell 1 at the left boundary on a uniform grid of spacing $h$ and consider the following second-order scheme:
\begin{eqnarray}
\frac{  u_1^{n+1} - u_1^{n}  }{\Delta t} +  a   \frac{  u_1^n -  u^{exact}( x_1 - h ,t^n)  }{h}  =  s( x_1, t^n ).
\label{special_scheme1}
\end{eqnarray}
This scheme retains an essential feature of the original scheme (\ref{original_scheme}). Namely, the truncation error is first order while the discretization error is second order, i.e., the truncation error is one order lower than the discretization error as is typical on irregular grids \cite{Boris_Jim_NIA2007-08,DiskinThomas:ANM2010}. It is designed to be independent of other cells, so that the discretization error can be easily derived. {\color{black} First, we substitute the exact solution into the scheme, 
\begin{eqnarray}
\frac{  u^{exact}( x_1 , t^{n+1}) - u^{exact}( x_1 ,t^n)  }{\Delta t} + a  \frac{u^{exact}( x_1  ,t^n) -  u^{exact}( x_1 - h ,t^n)  }{h}  =  s( x_1, t^n ) - {E}_1^n, 
\end{eqnarray}
where ${E}_1^n$ denotes the local truncation error at cell 1 at the time level $n$. Then, we subtract it from the scheme (\ref{special_scheme1}) to get 
\begin{eqnarray}
e_1^{n+1} = \left(  1 -  \mu  \right) e_1^{n}  +  \Delta t \, {E}_1^n, 
\label{error_eq}
\end{eqnarray}
where $\mu = a \Delta t / h$ is the CFL number, $e_1^{n} =u_1^{n} - u^{exact}( x_1 ,t^n) $ and $e_1^{n+1} =  u_1^{n+1}  -  u^{exact}( x_1 , t^{n+1})$ are the discretization errors at the time levels $n$ and $n+1$, respectively.} Here, we assume $\mu \le 1$ for convergence. For the sake of simplicity, we assume that the truncation error is independent of time and write
\begin{eqnarray}
{E}_1^n = C_1 ( \Delta t  + h), \quad  n=0, 1,2, \cdots, T_f / \Delta t,
\label{LTE_cell1}
\end{eqnarray}
which is derived straightforwardly by Taylor expanding Equation (\ref{special_scheme1}), where $C_1$ is a constant. Then, we obtain from Equation (\ref{error_eq}),
\begin{eqnarray}
e_1^{n} 
= \Delta t \,  C_1 ( \Delta t  + h)  \frac{ \displaystyle 1- \left(  1 - \mu  \right)^n    }{  \displaystyle 1-  \left(  1 - \mu   \right)  }
=
\frac{ \Delta t \,  C_1 ( \Delta t  + h)   }{ \mu } \left\{ 1- \left(  1 -  \mu   \right)^n  \right\} ,
\end{eqnarray}
which is the discretization error at the cell 1. 

In the case of accuracy verification with a single constant time step, we have $n=1$ and a constant $\Delta t$, and therefore, we find
\begin{eqnarray}
e_1^{1} 
= \Delta t  \,  C_1 ( \Delta t  + h) = O(h \Delta t) + O(\Delta t^2). 
\end{eqnarray}
{\color{black} It shows that the discretization error is proportional to $\Delta t$. Therefore, it can be made arbitrarily small by reducing $\Delta t$ but remains $O(h)$ in the grid refinement (i.e., as $h \rightarrow 0$). This explains the results shown in Figure \ref{fig:00_Li}, where extremely small first-order errors were obtained with an extremely small time step, $\Delta t = 10^{-8}$, at cells adjacent to a boundary on regular grids. This first-order error is originated from the first-order truncation error (\ref{LTE_cell1}). Note that the truncation
error is known to be one-order lower than the discretization error for
a scheme constructed in a linearly exact manner on an irregular stencil \cite{Boris_Jim_NIA2007-08,DiskinThomas:ANM2010,katz_sankaran:JSC_DOI,nishikawa_boundary_quadrature:JCP2015}, and the scheme here is such a scheme having an irregular stencil at the boundary cell. For irregular grids, therefore, first-order error convergence is expected everywhere on irregular grids, thus explaining the results for irregular grids found in Figures \ref{fig:00_L1} and \ref{fig:00_Li}. 
} 


On the other hand, if $n \ne 1$, we will have at the final time level $n = T_f/ \Delta t = T_f / ( \mu h / a )$,
\begin{eqnarray}
e_1^{n}  =
\frac{ \Delta t \,  C_1 ( \Delta t  + h)   }{ \mu }
\left\{ 1- \left(  1 - \mu   \right)^{\frac{a T_f}{\mu h}}  \right\}.
\end{eqnarray}
{\color{black} We can} expand it for $\mu \le 1$ as 
\begin{eqnarray}
e_1^{n} = 
\frac{ \Delta t \,  C_1 ( \Delta t  + h)   }{ \mu } 
\left[
1 - e^{ - \frac{a T_f}{h} }  +   \frac{a  \mu T_f }{2 h}e^{ - \frac{ a T_f}{h} }  +   O( \mu^2)
\right],
\end{eqnarray}
(see Ref.~\cite{Nishikawa_expansion:2023_RG_added} for the derivation of the Taylor series expansion of $ \left(  1 - \mu   \right)^{\frac{a T_f}{\mu h}} $). 
In the case of accuracy verification with $\Delta t = \mu h /a $ and a fixed $T_f$, where $\mu$ is a constant, we find
\begin{eqnarray}
e_1^{n} =  
\frac{  a \mu C_1 T_f }{2 a } \left( \frac{\mu}{a} +  1 \right) e^{ - \frac{a T_f}{h} } h 
+ 
 \frac{ C_1}{a} 
\left[
1 - e^{ - \frac{a T_f}{h} }
\right] \left( \frac{\mu}{a} +  1 \right)  h^2 
+ O( \mu^2).
\end{eqnarray}
Clearly, the first term introduces a first-order error if $a \, T_f / h$ is so small as to keep $ e^{ - \frac{a T_f}{h} }  = O(1)$. 
The second case mentioned in the previous section corresponds to $a=1$, $\mu=0.01$, $\Delta t =  \mu h $, $T_f = \mu h_c$, $h_c=1/8$, and $h= 1/ 2^{k+2}$, $k=1,2,3, \cdots, 6$, where
\begin{eqnarray}
\frac{ a T_f}{ h } = 0.00125 \times 2^{k+2} = \{ 0.01, 0.02, 0.04, 0.08, 0.16, 0.32 \},
\end{eqnarray}
which gives $ e^{ - \frac{ a T_f}{h} }  = O(1)$:
\begin{eqnarray}
e^{ - \frac{a T_f}{h} }  = \{ 0.99, 0.98, 0.96, 0.92, 0.85, 0.73 \}.
\label{expTfh_example_2}
\end{eqnarray}
This explains why we still observed first-order error convergence. 

It is emphasized again that these first-order errors appear because the truncation error is one order lower than the discretization error in the boundary stencil, where the grid is an irregular grid (i.e., stencils are different from one cell to another). In the interior of a regular grid, it is well known that the truncation and discretization errors have the same order \cite{Boris_Jim_NIA2007-08,DiskinThomas:ANM2010}, and therefore, we will observe second-order discretization errors in all cases. This explains why we obtained second-order error convergence in the $L_1$ norm as shown in Figure \ref{fig:00_L1}, i.e., when the second-order scheme in Section \ref{Fromm} has a second-order truncation error in the interior of the regular grids and the effect of first-order errors at boundary cells diminishes on average in the grid refinement. On irregular grids, the truncation error is known to be one order lower than the discretization errors for discretizations of a first-order hyperbolic partial differential equation (see, for example, Refs.~\cite{Boris_Jim_NIA2007-08,DiskinThomas:ANM2010,katz_sankaran:JSC_DOI,nishikawa_boundary_quadrature:JCP2015}). Therefore, the analysis implies that we would observe one-order lower error convergence generally for any high-order schemes. Indeed, we did observe second-order error convergence with a third-order convection scheme for irregular tetrahedral grids in three dimensions, when solving a time-dependent problem over a single tiny time step. 

In all cases, the first-order convergence may be understood as a coarse grid behavior. In the case of a constant $\Delta t$, the time domain is not refined and remains coarse. To observe a design order, we need to refine the time step along with the grid refinement. Then, we still have to make sure that the exponential factor in the first-order error will be made sufficiently small. This can be achieved by refining the grid, again implying that the grid is too coarse. However, there is a more practical approach that we will describe in the next section.

\vspace{-0.3cm}
\section{To Observe Design Orders}

\begin{figure}[h!]
  \begin{center}
      \begin{subfigure}[t]{0.32\textwidth}
        \includegraphics[width=\textwidth]{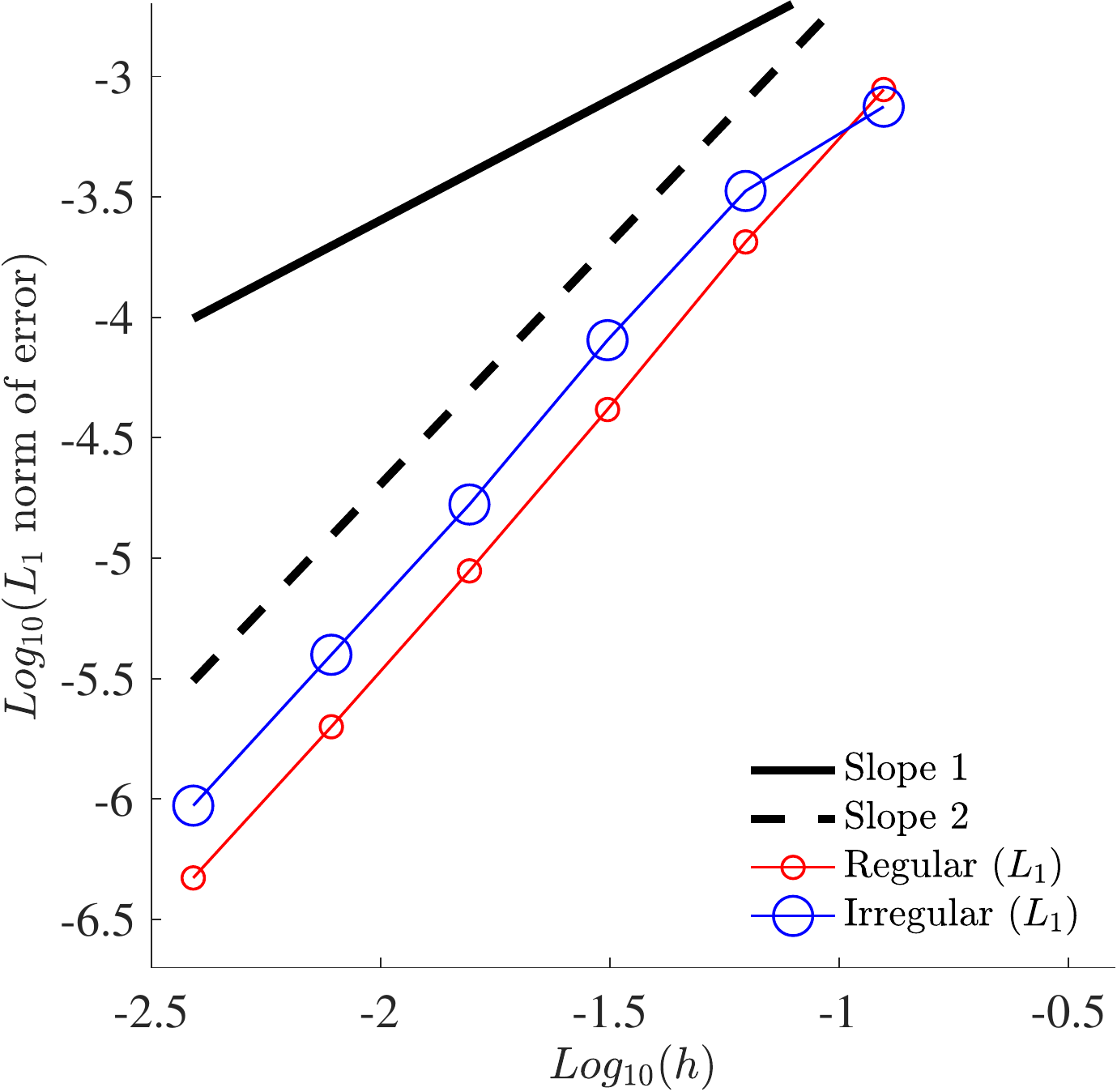}
          \caption{$L_1$ errors.}
          \label{fig:00_L1_largeDt}
      \end{subfigure}
      \begin{subfigure}[t]{0.32\textwidth}
        \includegraphics[width=\textwidth]{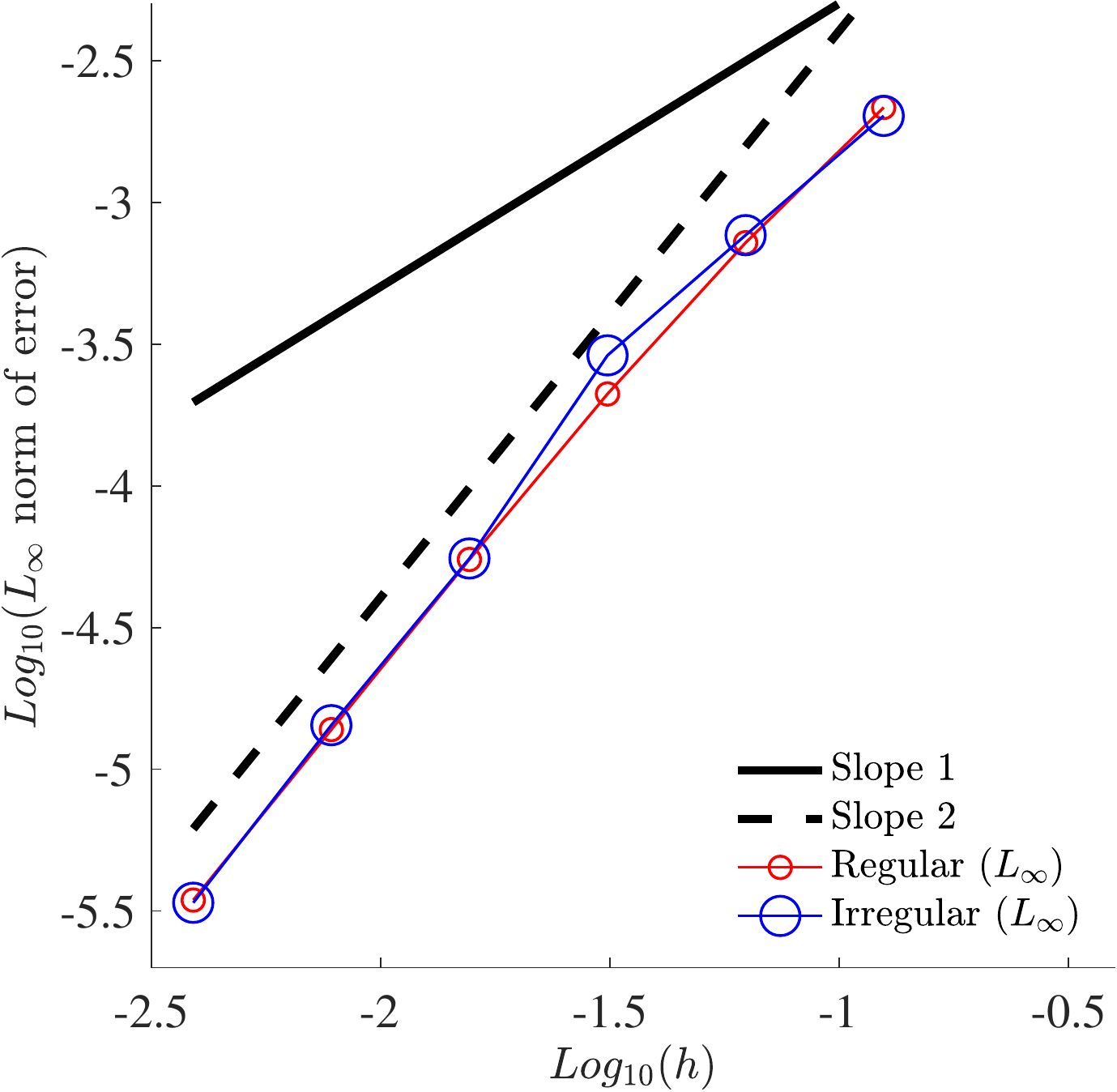}
          \caption{$L_\infty$ errors.}
          \label{fig:00_Li_largeDt}
      \end{subfigure} 
            \end{center}
            \caption{
\label{fig00_results}%
Accuracy verification results obtained with $\mu = 0.95$ and $T_f=\Delta t_c$, for a time-dependent problem on regular and irregular grids of $N = 2^{k+2}$ cells, $k=1,2, \cdots, 6$.
} 
\end{figure}

As the analysis has revealed, in order to observe design order of accuracy, we cannot just perform one constant time step, or we will encounter one-order lower discretization error convergence no matter how small the time step is. Instead, we should set $\Delta t = \mu h / a$ and choose the parameters, $T_f$ and $h$, such that the factor $e^{ - \frac{a T_f}{h} }$ is sufficiently small. Obviously, using finer grids (i.e., reducing $h$) is one possible solution, but this will require extremely fine grids. In the example corresponding to Equation (\ref{expTfh_example_2}), we have 
$e^{ - \frac{a T_f}{h} } = O(10^{-3})$ for $4,096$ cells and we will begin to observe second-order error convergence on this extremely fine grid. A more efficient approach would be to define the time step as
\begin{eqnarray}
\Delta t =  \frac{\mu h}{a},
\end{eqnarray}
with a constant $\mu$ of $O(1)$, and run the case toward $T_f = \Delta t_c =  \mu h_c / a$, where $h_c$ is the coarsest grid spacing. Then, we will have
\begin{eqnarray}
e^{ - \frac{a T_f}{h} }  = e^{ - n \mu   },
\end{eqnarray}
where $n = T_f / \Delta t = \Delta t_c /  \Delta t =  h_c / h$, which is, for our grids of spacing $h=1/2^{k+2}$ and $h_c = 1/8$, $n = 2^{k-1} = 1, 2, 4, 8, 16, 32$, resulting in
\begin{eqnarray}
e^{ - \frac{a T_f}{h} } = e^{ - n \mu   }  = \{ 0.37, \,   0.14, \,   1.8 \times 10^{-2},  \,  3.4 \times 10^{-4},   \,  1.1 \times 10^{-7},  \,  1.3 \times 10^{-14}  \},
\label{best_approach}
\end{eqnarray}
for $\mu  = 1$. Therefore, the first-order error is expected to vanish rapidly in the grid refinement and we would observe second-order error convergence using these grids. Numerical results obtained with $\mu = 0.95$ and $T_f = \Delta t_c $ are shown in Figures \ref{fig:00_L1_largeDt} and \ref{fig:00_Li_largeDt}, which confirm the prediction. In practice, the final time $T_f$ may be increased slightly, e.g., $T_f = 5 \Delta t_c $, so that the design order of accuracy can be observed right from the first two grid levels. Using this approach, we have successfully confirmed the design second- and third-order discretization error convergence for the three-dimensional cases mentioned earlier. 



\vspace{-0.3cm}
\section{Conclusions}

We have shown that accuracy verification can fail for time-dependent problems, exhibiting one order lower discretization error convergence no matter how small the time step is, or in other words, even if the time integration is nearly exact. As we have shown analytically, the discretization error near a boundary or in an irregular grid contains a lower-order term proportional to an exponential function of the ratio of the final time to a grid spacing, which will not vanish rapidly during the grid refinement unless the problem is set up carefully. It turns out that a correct and perhaps the most efficient approach to perform accuracy verification using a time-dependent problem is to integrate in time with the CFL number $\mu$ of $O(1) $ toward a final time defined by the time step on the coarsest grid level or longer by a factor of 5 or so, and perform the grid refinement from there. In this way, the exponential factor is forced to rapidly vanish, and a design order of accuracy is observed with a relatively small number of time steps on each grid level. 

 
\vspace{-0.2cm}
\section*{Acknowledgments}
 
The author gratefully acknowledges support by the Hypersonic Technology Project, through the Hypersonic Airbreathing Propulsion Branch of the NASA Langley Research Center. This work was funded under Contract No. 80LARC17C0004. The author would like to thank Boris Diskin (National Institute of Aerospace) for extensive discussions and valuable comments. The author would like to thank also Jeffery A. White (NASA Langley) for helpful comments.

\vspace{-0.2cm}
\bibliography{../../bibtex_nishikawa_database}
\bibliographystyle{unsrt}
\end{document}